\def\C{{\Cal C}}
\def\Z{{\Bbb Z}}
\def\E{{\Cal E}}

\documentstyle{amsppt}
\tolerance 3000
\pagewidth{5.5in}
\vsize7.0in
\magnification=\magstep1
\widestnumber \key{AAAAA}
\topmatter
\author Alex Iosevich, Nets Hawk Katz, and Terry Tao \endauthor
\address Department of Mathematics, Georgetown University, Washington,
DC 20057, USA \endaddress
\thanks Research supported in part by NSF grants.
\endthanks
\email iosevich\@math.georgetown.edu \endemail
\address Department of Mathematics, University of Illinois at Chicago,
Chicago, Illinois, 60607
\endaddress
\email nets\@math.uic.edu
\endemail
\address Department of Mathematics, University of California at Los Angeles,
Los Angeles, CA \endaddress
\email tao\@math.ucla.edu
\endemail
\title Convex bodies with a point of curvature do not have Fourier bases
\endtitle
\abstract We prove that no smooth symmetric convex body $\Omega$ with at
least
one point of non-vanishing Gaussian curvature can admit an orthogonal
basis of exponentials.  (The non-symmetric case was proven in \cite{Kol}).
This is further evidence of Fuglede's conjecture, which states that such a
basis is possible if and only if $\Omega$ can tile ${\Bbb R}^d$ by
translations.
\endabstract
\endtopmatter

\head Introduction and statement of results \endhead
\vskip.125in

Let $\Omega$ be a domain in ${\Bbb R}^d$, i.e., $\Omega$ is a Lebesgue
measurable subset of $\Bbb{R}^d$ with finite non-zero
Lebesgue measure. We say that a set $\Lambda \subset {\Bbb R}^d$ is
a {\it spectrum} of $\Omega$ if
${\{e^{2\pi i x \cdot \lambda} \}}_{\lambda \in \Lambda}$ is an orthogonal
basis of $L^2(\Omega)$.

\example{Conjecture A}\cite{Fug} A domain $\Omega$ admits a spectrum if and
only if it is
possible to tile ${\Bbb R}^d$ by a family of
translates of $\Omega$.
\endexample

Fuglede proved this conjecture under the additional assumption that the
tiling set or the spectrum are lattice subsets of ${\Bbb R}^d$. In general,
this conjecture is nowhere near resolution, even in
dimension one. It has been the subject of recent research,
see for example \cite{JoPe2} and \cite{LaWa}.

In this paper we shall address the following special case of Conjecture A.

\proclaim{Conjecture B} Suppose that $\Omega$ is a convex body with at least
one point of non-vanishing Gaussian curvature. Then $\Omega$ does not admit
a spectrum.
\endproclaim

The set $\Omega$ is called {\it symmetric} with respect to a point
$x_0 \in {\Bbb R}^d$ if $y \in \Omega$ implies that $2x_0-y \in \Omega$.
In \cite{Kol}, Kolountzakis proved Conjecture B under the assumption that
$\Omega$ is not symmetric with respect to any point.  However, the
case appears resistant to these methods.

In \cite{IKP}, (Theorem 1), the authors proved Conjecture B in the case
where $\Omega$ is the ball in ${\Bbb R}^d$, $d>1$. By generalizing the
arguments of \cite{IKP}, we show

\proclaim{Theorem 0.1} Suppose that $\Omega$ is a symmetric convex body in ${\Bbb R}^d$, 
$d \ge 2$. If the boundary of $\Omega$ is smooth, then $\Omega$ does not admit a spectrum.
The same conclusion holds in ${\Bbb R}^2$ if the boundary of $\Omega$ is piece-wise smooth, 
and has at least one point of non-vanishing Gaussian curvature. \endproclaim  

By Gauss-Bonnet theorem, a smooth hypersurface has at least one point of non-vanishing Gaussian 
curvature. Thus Conjecture B is true for smooth $\Omega$, and for piece-wise smooth 
$\Omega \subset {\Bbb R}^2$. 

The proof of Theorem 0.1 is based on the geometry of the set
$$Z_{\Omega}= \left\{\xi \in \Bbb R^d: \hat \chi_{\Omega}(\xi)=\int_\Omega
e^{-2\pi i \xi \cdot x}\ dx=0 \right\}. \tag0.1$$
The relevance of this set lies in the trivial observation that for any
spectrum $\Lambda$ of $\Omega$, we have
$$ \lambda - \lambda' \in Z_\Omega \hbox{ for all } \lambda, \lambda' \in
\Lambda, \lambda \neq \lambda'.  \tag0.2$$

For any $\eta \in {\Bbb R^d}$ and ball $B$, define the set
$X_{\Omega,\eta,B}$ by
$$ X_{\Omega,\eta,B} = Z_\Omega \cap B \cap (Z_\Omega - \eta) \cap (B -
\eta). \tag0.3$$

\definition{Definition} We say that a set $S \subset {\Bbb R}^d$ is 
$1$-separated if $\inf \{|x-y|: x,y \in S\}=1$. \enddefinition

Define the entropy $\E(X_{\Omega,\eta,B})$ to be the largest number of
1-separated points one can place inside $X_{\Omega,\eta,B}$.

Theorem 0.1 now follows immediately from the following two propositions.

\proclaim{Proposition 0.2} If $\Omega$ is a spectral set and $B$ is a ball
of radius $R \gg 1$, then there exists an $|\eta| \sim 1$ such that
$$ \E(X_{\Omega,\eta,B}) \sim R^d. \tag0.4$$
\endproclaim

\proclaim{Proposition 0.3} Let $\Omega$ be as in Theorem 0.1.  Then for
every $R \gg 1$ there exists a ball $B$ of radius $R$ and $\epsilon>0$ such that
$$ \E(X_{\Omega,\eta,B}) \lesssim R^{d-\epsilon} \tag0.5$$
for all $|\eta| \sim 1$.
\endproclaim

The proof of Proposition 0.3 will show that $\epsilon=1$ if the boundary of 
$\Omega$ is smooth. If $d=2$ and the boundary of $\Omega$ is piece-wise smooth, then
$\epsilon=\frac{1}{2}$. 

To illustrate these propositions we give two examples.  When $\Omega$ is a
cube, then $Z_\Omega$ is a union of hyperplanes, and $X_{\Omega,\eta,B}$
can be the union of $O(R)$ hyperplanes in $B$ with total entropy about
$R^d$.  However, when $\Omega$ is a sphere, $Z_\Omega$ is the union of
spheres, and $X_{\Omega,\eta,B}$ is the union of $O(R)$ $d-2$-dimensional
spheres in $B$, with total entropy about $R^{d-1}$.

Techniques similar to those used to prove Proposition 0.2 can be used to 
show the non-existence of spectra for other types of domains than convex 
bodies with a point of non-zero curvature. In a subsequent paper the authors 
will address this issue in the context of convex polygons.

\subhead Notation \endsubhead Throughout the paper, $a \sim b$, $a,b>0$,
means that there exist positive constants $c_1$ and $c_2$, such that $c_1a
\leq b \leq c_2a$. Similarly, $a \lesssim b$, $a,b>0$, means that there
exist a positive constant $c$, such that $a \leq cb$.

\head Proof of Proposition 0.2 \endhead
\vskip.125in

Let $\Omega$, $B$, $R$ be as in Proposition 0.2.  Let $\Lambda$ denote the 
putative spectrum for $\Omega$.

Let $B_1$, $B_2$ be any balls of radius $\sim R$ such that
$$B_1 - B_2 \subset B.\tag 1.1$$

Since $\hat \chi_\Omega$ is smooth and non-vanishing at the origin, we see
that $dist(0, Z_\Omega) \gtrsim 1$.  From (0.2) we thus have
$$|\lambda-\lambda'| \gtrsim 1 \hbox{ for all } \lambda, \lambda' \in
\Lambda,
\lambda \neq \lambda'. \tag1.2$$
We also have the density property
$$ \# B_i \gtrsim R^d \tag1.3$$
for $i=1,2$; see e.g. \cite{Lan}, \cite{Beu}, \cite{IoPe2}.

>From these two properties we may find $\lambda_1, \lambda_2 \in B_1$ so that
$$|\lambda_1-\lambda_2| \sim 1.\tag1.4$$
>From (0.2) and (1.1) we have
$$  \lambda - \lambda_1, \lambda - \lambda_2 \in Z_\Omega \cap B \tag1.5$$
for all $\lambda \in B_2 \cap \Lambda$.  We can re-arrange this as
$$ X_{\Omega,\lambda_2-\lambda_1,B} \supset (B_2 \cap \Lambda)-\lambda_2.
\tag1.6$$
The claim then follows from (1.2), (1.3), and (1.4).
\hfill$\square$

\head Proof of Proposition 0.3 \endhead
\vskip.125in

We first prove Proposition 0.3 in the case where the boundary of $\Omega$ is smooth. We shall 
then explain how the proof can be modified in two dimension to yield the conclusion of 
the theorem under the assumption that the boundary of $\Omega$ is piece-wise smooth and has
at least one point of non-vanishing Gaussian curvature. 

Let $\Omega$ be as in Theorem 0.1; we may assume that $\Omega$ is symmetric
around the origin.  Our main tool will be the method of stationary phase.

Our starting point is the formula
$$ \widehat{\chi}_{\Omega}(\xi)={(i|\xi|)}^{-1} \int_{\partial \Omega}
e^{2\pi i x \cdot \xi} \left( \frac{\xi}{|\xi|} \cdot n(x) \right)
d\sigma(x) \tag2.1$$
of Herz \cite{H}, where $n$ denotes the unit outward normal vector to
$\partial \Omega$ and $d\sigma$ denotes the Lebesgue measure on $\partial
\Omega$.

Let $x_0$ be a point on $\partial \Omega$ with non-vanishing Gaussian
curvature. By the symmetry assumption $-x_0$ is also in $\Omega$. Let $\psi$
be a smooth cutoff function supported in a small neighborhood of $x_0$. Let
$\C$ denote the cone of vectors normal to $\Omega$ on the support of $\psi$.
Integrating by parts we see that for $\xi \in \C$,
$$ (i|\xi|) \widehat{\chi}_{\Omega}(\xi)=
\int_{\partial \Omega} e^{2\pi i x \cdot \xi} \left( \frac{\xi}{|\xi|} \cdot
n(x) \right) (\psi(x)+\psi(-x))
d\sigma(x)+ O((1+|\xi|)^{-N}), \tag2.2$$
where $N$ is an arbitrary constant.

Fix $R \gg 1$, and let $B$ be a ball of radius $R$ in $\C$ which is a
distance $\sim R$ from the origin.  A stationary phase calculation (see e.g.
\cite{H}) gives
$$
\int_{\partial \Omega} e^{2\pi i x \cdot \xi} \frac{\xi}{|\xi|} \cdot n(x)
(\psi(x)+\psi(-x)) d\sigma(x) =
a(\xi) \cos \left(P(\xi) - \frac{\pi d}{4} \right) + O(R^{-\frac{d+1}{2}})
\tag2.3$$
for all $\xi \in B$, where $|a(\xi)| \sim R^{-\frac{d-1}{2}}$ and
$$ P(\xi) = \sup_{x \in \partial \Omega} x \cdot \xi. \tag2.4$$
Combining this with (2.2) we thus see that
$$ \left|\cos \left(P(\xi) - \frac{\pi d}{4} \right)\right| \lesssim R^{-1}
\tag 2.5$$
for all $\xi \in Z_\Omega \cap B$.

Fix $|\eta| \sim 1$. From (2.5) we see that
$$ dist(P(\xi+\eta) - P(\xi), \pi \Z) \lesssim R^{-1} \tag2.6$$
for all $\xi \in X_{\Omega,\eta,B}$.  From Taylor's theorem we thus have
$$ dist(\nabla P(\xi) \cdot \eta, \pi \Z) \lesssim R^{-1}. \tag2.7$$
Since $\nabla P(\xi) \cdot \eta$ is bounded, there are only a bounded number
of elements of $\pi \Z$ which are relevant in this distance calculation.
Since $\nabla P(\xi)$ is smooth, homogeneous of degree zero, and
non-degenerate on $B_R$ this means that $\xi/|\xi|$ lies within $O(R^{-1})$
of a finite number of $d-2$-dimensional surfaces in $S^{d-1}$ (which are
independent of $R$).  This implies the conclusion of Proposition 0.3 with 
$\epsilon=1$. 

If the boundary of $\Omega$ is piece-wise smooth, in any dimension, and has at least one point of non-vanishing Gaussian curvature, then $O(R^{-\frac{d+1}{2}})$ in $(2.3)$ must be replaced with $O(R^{-1})$. In two dimension this leads one to replace $O(R^{-1})$ in $(2.5)$, $(2.6)$, and $(2.7)$ by $O(R^{-\frac{1}{2}})$, which yields the conclusion of Proposition 0.3 with $\epsilon=\frac{1}{2}$. 

\hfill$\square$

\newpage
\vskip.25in
\head References \endhead

\ref \key Beu \by A. Beurling \paper Local harmonic analysis with some
applications to differential operators \yr 1966 \jour Some Recent
Advances
in the Basic Sciences, Academic Press \vol 1 \endref

\ref \key Fug \by B. Fuglede \paper Commuting self-adjoint partial
differential operators and a group theoretic problem \jour J. Funct.
Anal.
\yr 1974 \vol 16 \pages 101-121 \endref

\ref \key H \by C. S. Herz \paper Fourier transforms related to
convex sets \jour Ann. of Math. \vol 75 \yr 1962 \pages 81-92
\endref

\ref \key IKP \by A. Iosevich, N. Katz, and S. Pedersen
\paper Fourier basis and the Erd\H os distance problem
\jour Math. Research Letter \yr 1999 \vol 6 \pages \endref

\ref \key IoPe1 \by A. Iosevich and S. Pedersen \paper Spectral and
tiling
properties of the unit cube \jour Internat. Math Research Notices
\vol 16 \pages 819-828\yr 1998
\endref

\ref \key IoPe2 \by A. Iosevich and S. Pedersen \paper How large are
the spectral gaps?
\jour Pacific J. Math. (to appear) \vol \pages \yr 1998
\endref

\ref \key JoPe1 \by P. E. T. Jorgensen and S. Pedersen \paper Spectral
pairs in Cartesian coordinates
\jour J. Fourier Anal. Appl. (to appear) \vol \pages \yr 1998
\endref

\ref \key JoPe2 \by P. E. T. Jorgensen and S. Pedersen \paper Orthogonal
harmonic analysis of fractal measures
\jour ERA Amer. Math. Soc. \vol 4 \pages 35-42 \yr 1998
\endref

\ref \key Kol \by M. Kolountzakis \paper Non-symmetric convex domains
have no basis
of exponentials \jour (preprint) \vol \pages \yr 1999 \endref

\ref \key LaWa \by J. Lagarias and Y. Wang \paper
Spectral sets and factorizations of finite abelian groups
\jour J. Funct. Anal.
\vol 145 \pages 73-98 \yr 1997
\endref

\ref \key Lan \by H. Landau \paper Necessary density conditions for
sampling and interpolation of certain entire functions \jour Acta
Math.
\vol 117 \pages 37-52 \yr 1967 \endref
\enddocument